\documentclass[12pt]{article}
\usepackage[utf8]{inputenc}
\usepackage{amsmath,multicol}
\usepackage{amsthm}
\usepackage{amssymb}
\usepackage{graphicx}
\usepackage{color}
\usepackage{bm}
\usepackage{cases}

\allowdisplaybreaks[4]

\usepackage{mathrsfs}

\title{
Solutions
of quaternion-valued differential equations with or without commutativity
\thanks{Z. Cai was Supported by PhD research startup foundation of Hubei University of Technology (No.BSQD2019052).
K. I. Kou was supported by the Science and Technology Development Fund, Macau SAR (File no. FDCT/085/2018/A2) and University of Macau (File no. MYRG2019-00039-FST).
W. Zhang was supported by NSFC \# 11831012, \# 11821001 and \# 11771307.
}
}

\author{
{\normalsize
{\sc Zhenfeng Cai}\,$^a$,\quad {\sc Kit Ian Kou}\,$^b$
\thanks{Corresponding author: kikou@um.edu.mo
}
,\quad {\sc Weinian Zhang}\,$^c$
       }
\\
{\small  School of Science, Hubei University of Technology}
\\
{\small  Wuhan, Hubei 430068, China}
\\
{\small  Department of Mathematics, Faculty of Science and Technology} 
\\
{\small  University of Macau, Macau, China}
\\
{\small  School of Mathematics, Sichuan University}
\\
{\small  Chengdu, Sichuan 610064, China}
}

\date{}

\begin{document}
	
\maketitle

	\numberwithin{equation}{section}
	\newtheorem{lm}{Lemma}
	\newtheorem{thm}{Theorem}
	\newtheorem{remk}{Remark}
	\newtheorem{exam}{Example}
	\newtheorem{prop}{Proposition}
	\newtheorem{cor}{Corollary}
	\newtheorem{defn}{Definition}

\begin{abstract}
Most results on quaternion-valued differential equation (QDE) are based on
J. Campos and J. Mawhin's fundamental solution of exponential form for the homogeneous linear equation,
but their result requires a commutativity property.
In this paper we discuss with two problems:
What quaternion function satisfies the commutativity property?
Without the commutativity property, what can we do for the homogeneous equation?
We prove that the commutativity property actually requires quaternionic functions to be complex-like functions.
Without the commutativity property,
we reduce the initial value problem of the homogeneous equation
to a real nonautonomous nonlinear differential equation.

\vskip 0.1cm
		

\vskip 0.1cm

		
\end{abstract}

	\setlength\arraycolsep{4pt}
	\baselineskip 16pt
	\parskip 0.4cm


\section{Introduction}

The most fundamental form of quaternionic-valued differential equations is the following first order equation
\begin{eqnarray}
q'=F(t, q)
\label{qde}
\end{eqnarray}
with continuous $F: [0, T]\times \mathbb{H}\to \mathbb{H}$, where $\mathbb{H}$ denotes the set of all quaternions.
It corresponds, in terms of its components, to some systems of four first order real-valued differential equations.
The special interests in discussing quaternionic-valued differential equations not only
come from an analogue advantage of the complex formulation of planar vector fields, which enjoys
the multiplicative structure of complex numbers, but also arise from the
convenience of quaternions in formulating robotic manipulation (\cite{Ba01}), quaternion Frenet frames in differential geometry (\cite{Ha95}) and fluid mechanics (\cite{Gi02,Gi06}).
In 2006 Juan Campos and Jean Mawhin \cite{Cam-Mawhin} used topological degree methods to prove the existence of periodic
solutions for some special forms of (\ref{qde}).
Later Pawel Wilczynski (\cite{Wilcz09,Wilcz12}) gave some sufficient conditions for the existence of at least one or
two periodic solutions of the quaternion-valued Riccati equation and
detected solutions heteroclinic to the periodic ones.
Recently, optimization in quaternion-valued differential equations
was discussed in \cite{Xu16}.

In J. Campos and J. Mawhin's \cite{Cam-Mawhin}, the most-cited reference in papers of quaternion-valued differential equations,
the 1-dimensional linear QDE (abbreviation of quaternion-valued differential equation)
is of the form
\begin{eqnarray}
q'(t)= a(t) q(t) + f(t),
\label{lin_QDE}
\end{eqnarray}
where $a, f: \mathbb{R}\to \mathbb{H}$ are continuous on an interval $[0,T]$
and
$q: \mathbb{R}\to \mathbb{H}$ is unknown.
They noticed that the exponential of $q$, defined by
$$
\exp(q):=\sum_{n=0}^{\infty}\frac{q^n}{n!},
$$
satisfies the basic relation of exponential functions
$$
\exp(q)\exp(r) = \exp(q +r)
$$
if $r\in \mathbb{H}$ commutes with $q$, i.e., $rq=qr$,
and that
$
\{\exp(q(t))\}'=\{\exp(q(t))\} q'(t)
$
if $q$ is differentiable and commutes with its derivative, i.e.,
\begin{eqnarray}
q(t)q'(t) = q'(t)q(t)~~~\forall t\in [0,T].
\label{qq}
\end{eqnarray}
Then,
requiring the {\it commutativity property}
\begin{eqnarray}
a(t)A(t)=A(t)a(t),~~~\forall t\in [0,T],
\label{aA}
\end{eqnarray}
with $A(t):=A(0)+\int_0^t a(s)\,ds$ and a quaternionic constant $A(0)$,
in Proposition 3.1 and Corollary 3.3 of \cite{Cam-Mawhin}
they gave the fundamental solution $q(t)=\{\exp(A(t)-A(0))\}q(0)$
for the homogeneous equation
\begin{eqnarray}
q'(t)=a(t)q(t),
\label{h_lin_QDE}
\end{eqnarray}
which enables them to obtain
the Variation of Constant Formula
$$
q(t)=\exp(A(t)-A(0))\{q(0)+ \int_0^t \{\exp(A(0)-A(s))\}f(s)\,ds \}
$$
for the nonhomogeneous equation (\ref{lin_QDE})
in Corollary 3.4 of \cite{Cam-Mawhin}.
Based on the commutativity property, the necessary and sufficient conditions of the unique T-periodic solution for this equation were set up and an explicit formula in terms of a green function for this T-periodic solution was obtained in  \cite{Cam-Mawhin}. Furthermore, the fundamental matrix and solutions of two-dimensional QDE were given in \cite{Kou19} with a similar condition of (\ref{aA}).

The requirement of commutativity property (\ref{aA})
leads to two problems:
{\bf (Q1)} What quaternion function $a(t)$ satisfies (\ref{aA})?
{\bf (Q2)} Without (\ref{aA}), what can we do for the homogeneous equation (\ref{h_lin_QDE})?
In this paper we work on problems {\bf (Q1)} and {\bf (Q2)}.
For simple computation we only work with the special primitive function $A(t)$ of $a(t)$ such that $A(0)=0$.
In section 2 we give conditions for the quaternion function $a(t)$ under which the commutativity property (\ref{aA}) is true.
Moreover, we prove that (\ref{aA}) actually requires $a(t)$ to be a complex-like function.
In section 3, we discuss on the general solution of the homogeneous equation (\ref{h_lin_QDE}),
which is the most fundamental part in discussing differential equations with quaternion functions.
With the help of phase-angle representation of unit quaternions, the initial value problem of QDE (\ref{h_lin_QDE})
is turned to be a real nonautonomous nonlinear differential equation, called the decisive equation. Then,
an algorithm with Picard's iteration is applied to the decisive equation for solving QDE (\ref{h_lin_QDE}).
In section 4 we discuss in three special cases for
exact solutions
without using Picard's iteration.
We end this ppaer with examples in section 5.


\section{Commutativity property}

For the convenience of following discussion, we firstly introduce some symbols. We denote a quaternion $q$ by $q=q_0+q_1\mathbf{i}+q_2\mathbf{j}+q_3\mathbf{k}$, where $q_0,q_1,q_2,q_3$ are real numbers and $\mathbf{i},\mathbf{j},\mathbf{k}$ are the unit imaginary numbers which satisfy
\begin{equation*}
\mathbf{i}^2=\mathbf{j}^2=\mathbf{k}^2=\mathbf{ijk}=-1.
\end{equation*}

For the quaternion $h=h_0+h_1\mathbf{i}+h_2\mathbf{j}+h_3\mathbf{k}$, the inner product, norm are respectively defined by
\begin{equation*}
\begin{aligned}
<q,h>=q_0h_0+q_1h_1+q_2h_2+q_3h_3,
\\
|q|=\sqrt{<q,q>}=\sqrt{q_0^2+q_1^2+q_2^2+q_3^2}.
\end{aligned}
\end{equation*}
And the conjugate of $q$ is
\begin{equation*}
\overline{q}=q_0-q_1\mathbf{i}-q_2\mathbf{j}-q_3\mathbf{k},
\end{equation*}
and the inverse is
\begin{equation*}
q^{-1}=\frac{\overline{q}}{|q|^2}.
\end{equation*}

\begin{thm}\label{th1}
Let $a(t)=a_0(t) + a_1(t){\bf i} + a_2(t){\bf j} + a_3(t){\bf k}$ be continuous
and $B(t):=B(0)+\int_0^t a(s)\,ds$, where $B(0)=B_0(0)+B_1(0)\mathbf{i}+B_2(0)\mathbf{j}+B_3(0)\mathbf{k}$ is a certain quaternionic constant.
Then
$a(t)B(t)=B(t)a(t)$ for all $t\in [0,T]$
if and only if the real functions $a_1, a_2, a_3$ are proportional, i.e.,
\begin{eqnarray}
a_1(t):a_2(t):a_3(t)\equiv B_1(0):B_2(0):B_3(0),~~~\forall t\in [0,T].
\label{aaabbb}
\end{eqnarray}
Additionally,
$a(t)a'(t)=a'(t)a(t)$ for all $t\in [0,T]$
if $a$ is differentiable.
\label{comm1}
\end{thm}

Before proving the theorem, we need the following lemma.

\begin{lm}\label{lm1}
Let $p(t)=p_0(t)+\underline{p}(t)$ and $q(t)=q_0(t)+\underline{q}(t)$,
where $\underline{p}(t):=p_1(t)\mathbf{i}+p_2(t)\mathbf{j}+p_3(t)\mathbf{k}$ and
$\underline{q}(t):=q_1(t)\mathbf{i}+q_2(t)\mathbf{j}+q_3(t)\mathbf{k}$.
Then $p(t)q(t)=q(t)p(t)$ if and only if vector $\underline{p}(t)$ parallel to vector $\underline{q}(t)$
for each $t\in [0,T]$,
i.e. $\frac{p_1(t)}{q_1(t)}=\frac{p_2(t)}{q_2(t)}=\frac{p_3(t)}{q_3(t)}=\lambda(t)$ for all $t\in [0,T]$,
where $\lambda(t)$ is a certain real function.
\end{lm}

{\bf Proof.}
As defined in \cite[\S1.13, \S1.14]{MJ13},
\begin{equation*}
\begin{aligned}
p(t)q(t)
&=(p_0(t)+\underline{p}(t))(q_0(t)+\underline{q}(t))
\\
&=(p_0(t)q_0(t)-\underline{p}(t)\cdot\underline{q}(t))+p_0(t)\underline{q}(t)+q_0(t)\underline{p}(t)+\underline{p}(t)\times\underline{q}(t),
\end{aligned}
\end{equation*}
and
\begin{equation*}
\begin{aligned}
q(t)p(t)
&=(q_0(t)+\underline{q}(t))(p_0(t)+\underline{p}(t))
\\
&=(q_0(t)p_0(t)-\underline{q}(t)\cdot\underline{p}(t))+q_0(t)\underline{p}(t)+p_0(t)\underline{q}(t)+\underline{q}(t)\times\underline{p}(t).
\end{aligned}
\end{equation*}
where we know
\begin{equation*}
\begin{aligned}
&\underline{p}\cdot \underline{q}=p_1q_1+p_2q_2+p_3q_3, ~~~ \textrm{the dot product, and},
\\
&\underline{p}\times \underline{q}=\left|
                                                           \begin{array}{ccc}
                                                             \mathbf{i} & \mathbf{j} & \mathbf{k} \\
                                                             p_1 & p_2 & p_3 \\
                                                             q_1 & q_2 & q_3 \\
                                                           \end{array}
                                                         \right|, ~~~~~~~\,\textrm{the cross product},
\end{aligned}
\end{equation*}
and similarly we know $\underline{q}\cdot \underline{p}$ and $\underline{q}\times \underline{p}$.
Since $\underline{p}(t)\cdot\underline{q}(t)\equiv \underline{q}(t)\cdot\underline{p}(t)$, it follows that
$p(t)q(t)\equiv q(t)p(t)$
if and only if
$$
\underline{p}(t)\times\underline{q}(t)\equiv \underline{q}(t)\times\underline{p}(t),
$$
which means, by the counter-commutativity of the cross product, that
\begin{eqnarray}
\underline{p}(t)\times\underline{q}(t)\equiv 0.
\label{pxq}
\end{eqnarray}
It corresponds three trivial cases, which are: 1. $p_1(t)=p_2(t)=p_3(t)\equiv 0$ or $q_1(t)=q_2(t)=q_3(t)\equiv 0$; 2. $p_1(t)=q_1(t)\equiv 0$, $\frac{p_2(t)}{q_2(t)}=\frac{p_3(t)}{q_3(t)}$ or $p_2(t)=q_2(t)\equiv 0$, $\frac{p_1(t)}{q_1(t)}=\frac{p_3(t)}{q_3(t)}$ or $p_3(t)=q_3(t)\equiv 0$, $\frac{p_1(t)}{q_1(t)}=\frac{p_2(t)}{q_2(t)}$; 3. $p_1(t)=q_1(t)=p_2(t)=q_2(t)\equiv 0$,  or $p_1(t)=q_1(t)=p_3(t)=q_3(t)\equiv 0$, or $p_2(t)=q_2(t)=p_3(t)=q_3(t)\equiv 0$;  and a nontrivial case: where $p_i\not\equiv 0$ and $q_i\not\equiv0$, ($i=1,2,3$) and
$\underline{p}(t)\parallel \underline{q}(t)$,
i.e.,
\begin{eqnarray}
\frac{p_1(t)}{q_1(t)}=\frac{p_2(t)}{q_2(t)}=\frac{p_3(t)}{q_3(t)},
\label{p||q}
\end{eqnarray}
for all $t\in [0,T]\setminus D$, where $D$ is the set of zero points of $q_1(t),q_2(t)$ and $q_3(t)$. Since $\frac{0}{0}$ can be regarded as arbitrary number, the above cases can be generalized to all $t\in [0,T]$.
Let $\lambda(t)$ denote the common ratio of (\ref{p||q}), which is a real function. Thus the lemma is proved.

Having Lemma~\ref{lm1}, we are ready to prove the theorem.

{\bf Proof of Theorem~\ref{comm1}.}
Let
\begin{equation}\label{a-A}
A(t)
=\int_0^t a(s)\,ds:=A_0(t) + A_1(t){\bf i} + A_2(t){\bf j} + A_3(t){\bf k}.
\end{equation}
Then,
\begin{equation}\label{B(t)}
\begin{aligned}
B(t)&=B(0)+A(t):=B_0(t)+B_1(t)\mathbf{i}+B_2(t)\mathbf{j}+B_3(t)\mathbf{k}
\\
&=(B_0(0)+A_0(t)) + (B_1(0)+A_1(t)){\bf i} + (B_2(0)+A_2(t)){\bf j} + (B_3(0)+A_3(t)){\bf k},
\end{aligned}
\end{equation} and
\begin{eqnarray}
B_\ell'(t)=a_\ell(t),~~~\ell=0,\cdots3.
\label{AaA}
\end{eqnarray}
By Lemma \ref{lm1} or (\ref{pxq}), $a(t)B(t)\equiv B(t)a(t)$ if and only if
$$\left|
                                                           \begin{array}{ccc}
                                                             \mathbf{i} & \mathbf{j} & \mathbf{k} \\
                                                             a_1(t) & a_2(t) & a_3(t) \\
                                                             B_1(t) & B_2(t) & B_3(t) \\
                                                           \end{array}
                                                         \right|\equiv0
$$
for all $t\in [0,T]$.
For trivial case 1: $a_1(t)=a_2(t)=a_3(t)\equiv 0$ or $B_1(t)=B_2(t)=B_3(t)\equiv 0$ and trivial case 2: $a_1(t)=B_1(t)=a_2(t)=B_2(t)\equiv 0$  or $a_1(t)=B_1(t)=a_3(t)=B_3(t)\equiv 0$ or $a_2(t)=B_2(t)=a_3(t)=B_3(t)\equiv 0$; which naturally satisfy the result, i.e. $a_1(t),a_2(t)$ and $a_3(t)$ are proportional. For nontrivial case,
\begin{eqnarray}
\frac{B_1'(t)}{B_1(t)}= \frac{B_2'(t)}{B_2(t)}=\frac{B_3'(t)}{B_3(t)}=\frac{a_1(t)}{B_1(t)}
=\frac{a_2(t)}{B_2(t)}=\frac{a_3(t)}{B_3(t)},
\label{Diff-A}
\end{eqnarray}
for all $t\in [0,T],$
let $\lambda(t)$ be the common ratio of the above equalities. Since $a_i(t)=B_i'(t)$, any multiple root of $B_i(t)$ is still root of $a_i(t)$ and the difference of this two multiplicities is one. Hence $\lambda(t)$ is a piecewise continuous function which only has first order poles.
Solving real differential equations of (\ref{Diff-A}), we obtain
\begin{equation*}
B_1(t)=B_1(0) e^{\int_0^t\lambda(s)ds},~~
B_2(t)=B_2(0) e^{\int_0^t\lambda(s)ds},~~
B_3(t)=B_3(0) e^{\int_0^t\lambda(s)ds}.
\end{equation*}
By (\ref{AaA}),
\begin{equation*}
a_1(t)=B_1(0)\lambda(t)e^{\int_0^t\lambda(s)ds},~
a_2(t)=B_2(0)\lambda(t)e^{\int_0^t\lambda(s)ds},~
a_3(t)=B_3(0)\lambda(t)e^{\int_0^t\lambda(s)ds},
\end{equation*}
which means that $a_1(t):a_2(t):a_3(t)\equiv B_1(0):B_2(0):B_3(0)$
for all $t\in [0,T]$.

On the contrary, if (\ref{aaabbb}) holds, substituting (\ref{aaabbb}) in (\ref{a-A}),
we get $A_1(t):A_2(t):A_3(t)\equiv B_1(0):B_2(0):B_3(0)$.
Together with (\ref{B(t)}),
it follows that $B_1(t):B_2(t):B_3(t)=B_1(0):B_2(0):B_3(0)$,
which ensures
$$
\frac{a_1(t)}{B_1(t)}=\frac{a_2(t)}{B_2(t)}=\frac{a_3(t)}{B_3(t)}
$$
and implies $a(t)B(t)=B(t)a(t)$.

The additional resut is just a simple corollary with a substitution.
\qquad$\Box$


{\rm
Since for any continuous function $h(t)$, there exists a piecewise continuous function
$\lambda(t)={h(t)}/{\int_0^th(s)ds}$
such that $h(t)=\lambda(t)e^{\int_0^t\lambda(s)ds}$.  Condition (\ref{aaabbb}), i.e., $a_1(t):a_2(t):a_3(t)=B_1(0):B_2(0):B_3(0)$, is equivalent to the fact that
$a_1(t)=B_1(0)h(t)=B_1(0)\lambda(t)e^{\int_0^t\lambda(s)ds},a_2(t)=B_2(0)\lambda(t)e^{\int_0^t\lambda(s)ds}$ and $a_3(t)=B_3(0)\lambda(t)e^{\int_0^t\lambda(s)ds}$,
where $h(t)$ is a continuous function and $\lambda(t)$ is a piecewise continuous function.
Hence when $B_1^2(0)+B^2_2(0)+B_3^2(0)\neq0$,
$$
a(t)=a_0(t)+h(t)(B_1(0)\mathbf{i}+B_2(0)\mathbf{j}+B_3(0)\mathbf{k})=a_0(t)+g(t)\mathbf{I},
$$
where $g(t)=h(t)\sqrt{B_1^2(0)+B_2^2(0)+B_3^2(0)}$ is a continuous function and
\begin{eqnarray}
\mathbf{I}:=\frac{B_1(0)\mathbf{i}+B_2(0)\mathbf{j}+B_3(0)\mathbf{k}}{\sqrt{B_1^2(0)+B_2^2(0)+B_3^2(0)}}.
\label{III}
\end{eqnarray}

Thus,
the commutative condition (\ref{aA}) actually requires $a(t)$ to be a complex-like function.

\begin{thm}\label{co-complex}
The collection ${\cal I}:=\{w\in \mathbb{H}: w=a+b\mathbf{I}\}$, where
$\mathbf{I}$ is a certain pure quaternion defined by (\ref{III}), is a field.
Moreover, when $a(t)$ satisfies the commutative condition (\ref{aA}) i.e. $a(t)=a_0(t)+g(t)\mathbf{I}$,  and initial value $\mathbf{q}(0)$ belongs to ${\cal I}$, the Cauchy problem of equation (\ref{h_lin_QDE}) in quaternions is reduced to the Cauchy problem of the same equation in the field ${\cal I}$.
\label{Thm11}
\end{thm}

{\bf Proof.}
Since $\mathbf{I}^2=-1$, which means that $\mathbf{I}$ is a generalized image unit, we can easily prove that
if $ w_1:=a_1+b_1\mathbf{I}$ and $w_2:=a_2+b_2\mathbf{I}$ both belong to ${\cal I}$
then $w_1+w_2$, $w_1-w_2$, $w_1w_2$ and ${w_1}/{w_2}$ (as $w_2\neq0)$ belong to ${\cal I}$.
This implies that ${\cal I}$ is a field.

Furthermore, for each $ a+b\mathbf{I}\in {\cal I}$, there is a unique point $a+b\mathbf{i}\in\mathbb{C}$ and vice versa,
implying that ${\cal I}$ is isomorphic to $\mathbb{C}$.
When the Cauchy problem of equation (\ref{h_lin_QDE}) satisfies the commutative condition (\ref{aA})
and the initial value $\mathbf{q}(0)$ belongs to ${\cal I}$, the solution is
$$
\mathbf{q}(t)=\exp\left\{\left(\int_0^ta_0(s)\,ds\right)+\mathbf{I}\int_0^tg(s)\,ds\right\}\mathbf{q}(0),
$$
which belongs to ${\cal I}$.
Hence, this problem is reduced to the Cauchy problem of the same equation in the field ${\cal I}.$
\qquad$\Box$

\begin{remk}\label{rm1}
{\rm
If the Cauchy problem of equation (\ref{h_lin_QDE}) satisfies the commutative condition (\ref{aA})
but the initial value $\mathbf{q}(0)$ does not belong to ${\cal I}$,
then the solution does not belong to field ${\cal I}$.
In this case, however,
we can still find the general solution
from the Cauchy problem of the same equation in the field ${\cal I}$,
 which is $\mathbf{q}(t)=\exp\{(\int_0^ta_0(s)\,ds)+\mathbf{I}\int_0^tg(s)\,ds\}\mathbf{q}(0)$, where $\exp\{(\int_0^ta_0(s)\,ds)+\mathbf{I}\int_0^tg(s)\,ds\}$ belongs to ${\cal I}$.
}
\end{remk}

For example,
consider the homogeneous linear equation (\ref{h_lin_QDE}) associated with the initial condition
$q(0)=\mathbf{i}$, where $a(t)=t^2+t(\mathbf{i}+2\mathbf{j}+3\mathbf{k})$.
Clearly,
$a_1(t):a_2(t):a_3(t)=1:2:3$. By Theorem \ref{th1}, the commutative condition (\ref{aA}) holds.
According to Theorem \ref{co-complex}, $a(t)=t^2+\sqrt{14}t\frac{\mathbf{i}+2\mathbf{j}+3\mathbf{k}}{\sqrt{14}}=t^2+\sqrt{14}t\mathbf{I}$, and the solution
$$
q(t)=e^{\int_0^t a(s)ds}q(0)=e^{\int_0^t a(s)ds}\mathbf{i}=e^{\frac{t^3}{3}+\frac{\sqrt{14}t^2}{2}\mathbf{I}}\mathbf{i},
$$
where $\mathbf{I}:=\frac{\mathbf{i}+2\mathbf{j}+3\mathbf{k}}{\sqrt{14}}$. Indeed,
$$
q(t)=e^{\frac{t^3}{3}}\bigg\{-\frac{1}{\sqrt{14}}\sin \left(\frac{\sqrt{14}}{2}t^2\right)+\cos \left(\frac{\sqrt{14}}{2}t^2 \right)\mathbf{i}+\frac{3}{\sqrt{14}}\sin \left(\frac{\sqrt{14}}{2}t^2 \right)\mathbf{j}-\frac{2}{\sqrt{14}}\sin \left(\frac{\sqrt{14}}{2}t^2 \right)\mathbf{k}\bigg\},
$$
and it does not belong to the field ${\cal I}:=\left\{w\in \mathbb{H}: w=a+b\frac{\mathbf{i}+2\mathbf{j}+3\mathbf{k}}{\sqrt{14}}\right\}$.


\section{
Solutions without commutativity}

In this section we discuss equation (\ref{h_lin_QDE}) without the commutative condition (\ref{aA}).
In order to simplify the discussion, we note the following.

\begin{lm}
\label{lm2-thm}
Let $a(t):=a_0(t)+\underline{a}(t)$, where $\underline{a}(t):=a_1(t)\mathbf{i}+a_2(t)\mathbf{j}+a_3(t)\mathbf{k}$.
If $y_1(t)$ is a real-valued solution of the equation $y_1'(t)=a_0(t)y_1(t)$
and $y_2(t)$ is a solution of QDE $y_2'(t)=\underline{a}(t)y_2(t)$.
Then the function $Y(t):=y_1(t)y_2(t)$ is a solution of QDE (\ref{h_lin_QDE}).
\end{lm}

{\bf Proof.}
As supposed, $y_1'(t)=a_0(t)y_1(t)$ and $y_2'(t)=\underline{a}(t)y_2(t)$.
Then
$$
Y'(t)
=y'_1(t)y_2(t)+y_1(t)y'_2(t)
=a_0(t)y_1(t)y_2(t)+y_1(t)\underline{a}(t)y_2(t).
$$
Being real valued, $y_1(t)$ is commutative with $\underline{a}(t)$.
It follows that
$$
Y'(t)=\big(a_0(t)+\underline{a}(t)\big)y_1(t)y_2(t)=a(t)Y(t).
$$
This implies that $Y(t)$ is a solution of the equation $q'(t)=a(t)q(t)$.
\qquad$\Box$

Since real solutions of the equation $y_1'(t)=a_0(t)y_1(t)$ can be found easily,
Lemma \ref{lm2-thm} shows that the key of solving QDE (\ref{h_lin_QDE}) is to find solutions of the equation
\begin{equation}
\label{PureQDE}
y'(t)=\underline{a}(t)y(t).
\end{equation}
Thus, we can focus on this QDE with $a_0(t)=0$, i.e., $a(t)=\underline{a}(t)=a_1(t)\mathbf{i}+a_2(t)\mathbf{j}+a_3(t)\mathbf{k}$.

\begin{lm}
\label{lm22-thm}
If
$y(t):=q_0(t)+q_1(t)\mathbf{i}+q_2(t)\mathbf{j}+q_3(t)\mathbf{k}$
satisfies equation (\ref{PureQDE}) with $a_0(t)=0$, then
$|y(t)|$ is a constant for all $t$.
\end{lm}

{\bf Proof.}
Let
$q(t)=q_0(t)+q_1(t)\mathbf{i}+q_2(t)\mathbf{j}+q_3(t)\mathbf{k}$ and
$a(t)=a_0(t)+a_1(t)\mathbf{i}+a_2(t)\mathbf{j}+a_3(t)\mathbf{k}$.
Then we can generally reduce QDE (\ref{h_lin_QDE}) to the ODE form
\begin{equation*}\label{ODE-1}
\left\{\begin{aligned}
&q'_0=a_0(t)q_0 -a_1(t)q_1-a_2(t)q_2-a_3(t)q_3,
\\
&q'_1=a_1(t)q_0+a_0(t)q_1 -a_3(t)q_2+a_2(t)q_3,
\\
&q'_2=a_2(t)q_0+a_3(t)q_1+a_0(t)q_2 -a_1(t)q_3,
\\
&q'_3=a_3(t)q_0-a_2(t)q_1+a_1(t)q_2+a_0(t)q_3,
\end{aligned}
\right.
\end{equation*}
or equivalently the form
\begin{eqnarray}
\vec{q}\;'= M(t) \vec{q},
\label{real-vector-equation}
\end{eqnarray}
where $\vec{q}:=(q_0,q_1,q_2,q_3)^T$ is a 4-dimensional real vector and
\begin{equation*}
M(t)
:=\left(
\begin{array}{rrrr}
a_0(t) & -a_1(t) & -a_2(t) & -a_3(t) \\
a_1(t) & a_0(t) & -a_3(t) & a_2(t) \\
a_2(t) & a_3(t) & a_0(t) & -a_1(t) \\
a_3(t) & -a_2(t) & a_1(t) & a_0(t) \\
\end{array}
\right)
\end{equation*}
is a $4\times 4$ real matrix.
Hence,
for QDE (\ref{PureQDE})
with $a_0(t)=0$,
we only need to consider the 4-dimensional first order ODE with variable coefficients
\begin{equation}
\label{ODE-2}
\left\{\begin{aligned}
&q'_0=\quad \qquad-a_1(t)q_1-a_2(t)q_2-a_3(t)q_3,
\\
&q'_1=a_1(t)q_0\qquad\qquad -a_3(t)q_2+a_2(t)q_3,
\\
&q'_2=a_2(t)q_0+a_3(t)q_1 \qquad\qquad-a_1(t)q_3,
\\
&q'_3=a_3(t)q_0-a_2(t)q_1+a_1(t)q_2\qquad\qquad,
\end{aligned}
\right.
\end{equation}
a special case of equation (\ref{real-vector-equation}).
Using equations of (\ref{ODE-2}) to replace
$q_0',q_1',q_2'$ and $q_3'$, we compute
\begin{eqnarray*}
(q_0^2+q_1^2+q_2^2+q_3^2)'
&=&2q_0q_0'+2q_1q_1'+2q_2q_2'+2q_3q_3'
\\
&=&
2q_0\big(-a_1(t)q_1-a_2(t)q_2-a_3(t)q_3\big)
\\
&&\quad
+2q_1\big(a_1(t)q_0-a_3(t)q_2+a_2(t)q_3\big)
\\
&&\quad
+2q_2\big(a_2(t)q_0+a_3(t)q_1-a_1(t)q_3\big)
\\
&&\quad
+2q_3\big(a_3(t)q_0-a_2(t)q_1+a_1(t)q_2\big)
\\
&=&
0,
\end{eqnarray*}
which implies that $|y(t)|=\sqrt{q_0^2+q_1^2+q_2^2+q_3^2}=C_0$, a non-negative real constant. This completes the proof.
\qquad$\Box$

As shown in \cite{Kou19}, the general linear QDE has the same superposition principle, i.e.,
any linear combination of solutions is still a solution, as real and complex linear differential equations
although
those solutions do not form a linear space but a right module over the divisible quaternionic ring.
Being a special case, QDE (\ref{PureQDE}) has the
general solution $y(t)=q(t)\mathbf{C}$,
where $q(t)$ is a nontrivial solution of QDE (\ref{PureQDE}) and $\mathbf{C}$ is an arbitrary quaternionic constant.
Since the solution $q(t)$ of QDE (\ref{PureQDE}) with $a_0(t)\equiv 0$
has a constant norm by Lemma \ref{lm22-thm}, equation (\ref{PureQDE})
has a nontrivial solution $y(t)$ such that
 $|y(t)|=1$, which corresponds to $\mathbf{C}:=|q(t)|^{-1}$.
Furthermore, for the initial value problem
\begin{equation}
\label{PureQDE-IVP}
y'(t)=\underline{a}(t)y(t),\quad y(t_0)=\mathbf{C}_0,
\end{equation}
the solution can be expressed by $y(t)\mathbf{C}_0$, where $y(t)$ is not only the solution of  QDE (\ref{PureQDE}) with unit norm, but also has the unit initial value which is $y(t_0)=1$.
Hence the key of solving QDE (\ref{PureQDE-IVP}) is turned to solve the initial problem
\begin{equation}
\label{PureQDE-IVP-1}
y'(t)=\underline{a}(t)y(t),\quad y(t_0)=1.
\end{equation}
For this purpose, we need the following lemma.

\begin{lm}
{\bf (Corollary 2 of \cite{Bulow2001})}    
\label{Exponential-form}
Every unit quaternion $q=q_0+q_1\mathbf{i}+q_2\mathbf{j}+q_3\mathbf{k}$  can be represented as $q=e^{\mathbf{i}\theta_1}e^{\mathbf{j}\theta_2}e^{\mathbf{k}\theta_3}$,
where $(\theta_1,\theta_2, \theta_3)
\in(-\pi,\pi]\times[-{\pi}/{4},{\pi}/{4}]\times(-{\pi}/{2},{\pi}/{2}]$, being the phase of $q$, is determined by either
\begin{equation*}
\theta_2=\frac{1}{2}\arcsin(2(q_0q_2+q_1q_3))\in (-\pi/4, \pi/4)
\end{equation*}
and
\begin{equation}\label{theta-1}
\left\{
\begin{aligned}
\theta_1&=\mathrm{atan}2 \left(\frac{q_0+q_1+q_2+q_3}{q_0-q_1+q_2-q_3}\right)-\frac{\pi}{4}-\theta_3,
\\
\theta_3&=\frac{1}{2}\mathrm{atan}2(2(q_0q_3-q_1q_2),q_0^2+q_1^2-q_2^2-q_3^2),
\end{aligned}
\right.
\end{equation}
or
\begin{equation}\label{theta-2}
\theta_2=\frac{\pi}{4} \quad and \quad
\theta_1=\mathrm{atan}2 \left(\frac{q_0+q_1+q_2+q_3}{q_0-q_1+q_2-q_3} \right)-\frac{\pi}{4}-\theta_3,
\end{equation}
or
\begin{equation}\label{theta-3}
\theta_2=-\frac{\pi}{4} \quad and \quad
\theta_1=\mathrm{atan}2 \left(\frac{q_0+q_1-q_2-q_3}{q_0-q_1-q_2+q_3} \right)-\frac{\pi}{4}+\theta_3,
\end{equation}
with arbitrary $\theta_3\in (-{\pi}/{2},{\pi}/{2}]$. \\

%
%
\end{lm}
\begin{remk}
{\rm
{\bf (i)}
Due to the periodicity of  $e^{\mathbf{i}\theta_1}$, one can find a proper $\theta_1$ belongs to $(-\pi,\pi]$  by adding or subtracting $2\pi$ in equations (\ref{theta-1}), (\ref{theta-2}) and (\ref{theta-3}).
{\bf (ii)}
The function $\mathrm{atan2}$ is defined as: $\mathrm{atan2}(b,a):=\arctan(b/a)$ for $a>0$,
$:= \pm \frac{\pi}{2}$ for $a=0$, $:=\arctan(b/a)-\pi$ for $a<0$ and $b<0$, and $:=\arctan(b/a)+\pi$ for $a<0$ and $b\geq0$.
{\bf (iii)}
The phrase ``almost uniquely'' is used since there are two singular cases $\theta_2=\pm\frac{\pi}{4}$.
For $\theta_2\in (-\pi/4,\pi/4)$ those formulae in (\ref{theta-1})
give a unique $(\theta_1,\theta_2,\theta_3)$
but
for $\theta_2=\pm\pi/4$ there are infinitely many possiblities of $(\theta_1,\theta_2,\theta_3)$
because
all values of $\theta_1$ and $\theta_3$ satisfying $\theta_1\pm\theta_3=C$ for some constant $C$ fulfill this expression of $q$.
In these singular cases,
as shown in \cite{BulowPHD,Bulow2001},
engineering scientists would like to set $\theta_3=0$,
which implies that $q$ has expressions
$$
\theta_1=\mathrm{atan}2 \left(\frac{q_0+q_1+q_2+q_3}{q_0-q_1+q_2-q_3} \right)-\frac{\pi}{4} \mbox{ and } \theta_1=\mathrm{atan}2 \left(\frac{q_0+q_1-q_2-q_3}{q_0-q_1-q_2+q_3}\right)-\frac{\pi}{4}
$$
uniquely for $\theta_2=\frac{\pi}{4}$ and $\theta_2=-\frac{\pi}{4}$ respectively.
}
\end{remk}

\begin{remk}\label{rmketc}
{\rm
Due to the non-commutativity of quaternions, we have $(\mathbf{i}\theta_1)(\mathbf{j}\theta_2)\neq (\mathbf{j}\theta_2)(\mathbf{i}\theta_1)$ and $e^{\mathbf{i}\theta_1}e^{\mathbf{k}\theta_2}e^{\mathbf{j}\theta_3}\neq e^{\mathbf{i}\theta_1}e^{\mathbf{j}\theta_2}e^{\mathbf{k}\theta_3}$.
We similarly see that every unit quaternion $q$ has another five expressions, i.e.,
$q=e^{\mathbf{i}\theta_1}e^{\mathbf{j}\theta_2}e^{\mathbf{k}\theta_3}$,
$q=e^{\mathbf{j}\theta_1}e^{\mathbf{i}\theta_2}e^{\mathbf{k}\theta_3}$,
$q=e^{\mathbf{j}\theta_1}e^{\mathbf{k}\theta_2}e^{\mathbf{i}\theta_3}$,
$q=e^{\mathbf{k}\theta_1}e^{\mathbf{i}\theta_2}e^{\mathbf{j}\theta_3}$, and
$q=e^{\mathbf{k}\theta_1}e^{\mathbf{j}\theta_2}e^{\mathbf{i}\theta_3}$.
}
\end{remk}

Now we return to find a solution with the unit norm for the initial value problem (\ref{PureQDE-IVP-1}).
Since every unit quaternionic function $q(t)$ satisfies
$q^{-1}(t)={\overline{q}(t)}/{|q(t)|^2}=\overline{q}(t)$
as a complex unit one,
}
we see that
 solutions of equation (\ref{PureQDE}) with unit norms satisfy
 \begin{equation}\label{PureQDE-2}
q'(t)\overline{q}(t)=\underline{a}(t),\quad q(t_0)=1.
\end{equation}
By Lemma \ref{Exponential-form}, we only need to find a solution of
(\ref{PureQDE-2}) in the form $e^{\mathbf{i}\theta_1(t)}e^{\mathbf{j}\theta_2(t)}e^{\mathbf{k}\theta_3(t)}$.

\begin{thm}\label{Corr-ODE}
For $\theta_2\neq\pm\frac{\pi}{4}$,
the quaternion
$q(t)=e^{\mathbf{i}\theta_1(t)}e^{\mathbf{j}\theta_2(t)}e^{\mathbf{k}\theta_3(t)}$ is a nontrivial solution of QDE (\ref{PureQDE-IVP-1})
with the unit initial value
if the vector
$(\theta_1(t),\theta_2(t),\theta_3(t))^T$ is a real solution of the system
\begin{equation}\label{Corrspond-ODE}
\left\{\begin{aligned}
\theta'_1(t)&=a_1(t)+\sin(2\theta_1)\tan(2\theta_2)a_2(t)-\cos(2\theta_1)\tan(2\theta_2)a_3(t),
\\
\theta'_2(t)&=\cos(2\theta_1)a_2(t)+\sin(2\theta_1)a_3(t),
\\
\theta'_3(t)&=-\frac{\sin(2\theta_1)}{\cos(2\theta_2)}a_2(t)+\frac{\cos(2\theta_1)}{\cos(2\theta_2)}a_3(t),
\end{aligned}\right.
\end{equation}
with $(\theta_1(t_0),\theta_2(t_0),\theta_3(t_0))=(0,0,0)$.
\end{thm}

{\bf Proof.}
Using the form $q(t)=e^{\mathbf{i}\theta_1(t)}e^{\mathbf{j}\theta_2(t)}e^{\mathbf{k}\theta_3(t)}$, one can compute
\begin{eqnarray*}
q'(t)
&=&(e^{\mathbf{i}\theta_1(t)}e^{\mathbf{j}\theta_2(t)}e^{\mathbf{k}\theta_3(t)})'
\\
&=&\mathbf{i}\theta'_1(t)e^{\mathbf{i}\theta_1(t)}e^{\mathbf{j}\theta_2(t)}
e^{\mathbf{k}\theta_3(t)}+e^{\mathbf{i}\theta_1(t)}\mathbf{j}\theta'_2(t)e^{\mathbf{j}\theta_2(t)}e^{\mathbf{k}\theta_3(t)}
+e^{\mathbf{i}\theta_1(t)}e^{\mathbf{j}\theta_2(t)}\mathbf{k}\theta_3'(t)e^{\mathbf{k}\theta_3(t)},
\\
\overline{q}(t)
&=&
e^{-\mathbf{k}\theta_3(t)}e^{-\mathbf{j}\theta_2(t)}e^{-\mathbf{i}\theta_1(t)},
\end{eqnarray*}
from which we get
\begin{align}
q'(t)\overline{q}(t)&=\big(\mathbf{i}\theta'_1(t)
e^{\mathbf{i}\theta_1(t)}e^{\mathbf{j}\theta_2(t)}e^{\mathbf{k}\theta_3(t)}+e^{\mathbf{i}\theta_1(t)}
\mathbf{j}\theta'_2(t)e^{\mathbf{j}\theta_2(t)}e^{\mathbf{k}\theta_3(t)}+e^{\mathbf{i}\theta_1(t)}
e^{\mathbf{j}\theta_2(t)}\mathbf{k}\theta_3'(t)e^{\mathbf{k}\theta_3(t)}\big)\notag
\\
&\quad\big(e^{-\mathbf{k}\theta_3(t)}e^{-\mathbf{j}\theta_2(t)}e^{-\mathbf{i}\theta_1(t)}\big)\notag
\\
&=\mathbf{i}\theta'_1(t)+e^{\mathbf{i}\theta_1(t)}\mathbf{j}\theta'_2(t)e^{-\mathbf{i}\theta_1(t)}
+e^{\mathbf{i}\theta_1(t)}e^{\mathbf{j}\theta_2(t)}\mathbf{k}\theta_3'(t)e^{-\mathbf{j}\theta_2(t)}
e^{-\mathbf{i}\theta_1(t)}\label{ijk-1}
\\
&=\bigg(\theta_1'(t)+\theta_3'(t)\sin(2\theta_2(t))\bigg)\mathbf{i}
+\bigg(\theta_2'(t)\cos(2\theta_1(t))-\theta'_3(t)\sin(2\theta_1(t))\cos(2\theta_2(t))\bigg)\mathbf{j}
\notag\\
&\ \
+\bigg(\theta_2'(t)\sin(2\theta_1(t))+\theta'_3(t)\cos(2\theta_1(t)) \cos(2\theta_2(t))\bigg)\mathbf{k}\label{ijk-2}.
\end{align}
Thus, we rewrite equation (\ref{PureQDE-2}) as the equivalent real system
\begin{equation}\label{important-form}
\left\{\begin{aligned}
&\theta'_1(t)+\theta'_3(t)\sin(2\theta_2(t))=a_1(t),\qquad\qquad\qquad\qquad
\\
&\theta_2'(t)\cos(2\theta_1(t))-\theta'_3(t)\sin(2\theta_1(t))\cos(2\theta_2(t))=a_2(t),
\\
&\theta_2'(t)\sin(2\theta_1(t))+\theta'_3(t)\cos(2\theta_1(t))\cos(2\theta_2(t))=a_3(t).
\end{aligned}
\right.
\end{equation}
By Lemma \ref{Exponential-form}, the initial value $q(t_0)=1$ corresponds to $(\theta_1(t_0),\theta_2(t_0),\theta_3(t_0))=(0,0,0)$.
Furthermore,
regarding (\ref{important-form}) as a linear system of $\theta_1'(t),\theta'_2(t)$ and $\theta_3'(t)$,
we obtain its coefficients matrix
\begin{equation*}
A:=\left(\begin{array}{cccc}
    1 &    0               & \quad\qquad\sin(2\theta_2) \\
    0 &    \cos(2\theta_1)   & -\sin(2\theta_1)\cos(2\theta_2)\\
    0 &    \sin(2\theta_1)   & \cos(2\theta_1)\cos(2\theta_2)
\end{array}\right).
\end{equation*}
Since $\theta_2(t)\neq\pm\frac{\pi}{4}$,  which corresponds to the invertible case of $A$,
we can compute $\det(A)=\cos(2\theta_2(t))$ and thus
\begin{equation*}
A^{-1}=\left(\begin{array}{cccc}
    1 &    \sin(2\theta_1)\tan(2\theta_2)              & \quad\qquad-\cos(2\theta_1)\tan(2\theta_2) \\
    0 &    \cos(2\theta_1)   & -\sin(2\theta_1)\\
    0 &    -\frac{\sin(2\theta_1)}{\cos(2\theta_2)}   & \frac{\cos(2\theta_1)}{\cos(2\theta_2)}
\end{array}\right),
\end{equation*}
which reduces (\ref{important-form}) to (\ref{Corrspond-ODE}).
The proof is completed.
\qquad$\Box$


\begin{remk}\label{singular-case}
{\rm
Theorem \ref{Corr-ODE} is effective only in the region without $\theta_2(t)=\pm\frac{\pi}{4}$ because $A$ is not invertible there.
}
\end{remk}

By Theorem \ref{Corr-ODE}, the key of solving QDE (\ref{PureQDE})
is to find solutions of equation (\ref{Corrspond-ODE}),
called the {\it decisive equation},
which can be expressed as
\begin{equation}\label{Corrspond-ODE2}
\frac{d\mathbf{x}}{dt}=\mathbf{f}(t,\mathbf{x}),
\end{equation}
 where
\begin{equation*}
\mathbf{x}=\left(
  \begin{array}{c}
    \theta_1(t) \\
    \theta_2(t) \\
    \theta_3(t) \\
  \end{array}
\right),
\mathbf{f}(t,\mathbf{x})=\left(
  \begin{array}{r}
    a_1(t)+\sin(2\theta_1)\tan(2\theta_2)a_2(t)-\cos(2\theta_1)\tan(2\theta_2)a_3(t) \\
    \cos(2\theta_1)a_2(t)+\sin(2\theta_1)a_3(t) \\
    -\frac{\sin(2\theta_1)}{\cos(2\theta_2)}a_2(t)+\frac{\cos(2\theta_1)}{\cos(2\theta_2)}a_3(t) \\
  \end{array}
\right).
\end{equation*}
The existence and uniqueness of solutions of this equation
are given in the following.

\begin{thm}\label{existence}
Let $a_1(t),a_2(t)$ and $a_3(t)$ be continuous functions on the interval $[t_0-a, t_0+a]$.
the decisive equation (\ref{Corrspond-ODE2}) with the initial value
$\mathbf{x}(t_0)=(0,0,0)$
has a unique solution on the interval $I=[t_0-h,t_0+h]$,
where $h:=\min\{a,\frac{b}{M}\}$, $M:=\max\{|\mathbf{f}(t,\mathbf{x})|:|t-t_0|\leq a, |\mathbf{x}|\leq b\}$ and $b$ is a positive constant less than $\pi/4$,.
\end{thm}

{\bf Proof.}
Let
$D:=\{(t,\mathbf{x})||t-t_0|\leq a,|\mathbf{x}(t)|\leq b\}$,
which is a closed domain in $\mathbb{R}\times (-\pi,\pi]\times [-\frac{\pi}{4},\frac{\pi}{4}]\times (-\frac{\pi}{2},\frac{\pi}{2}]$  without the points with $\theta_2(t)=\pm \frac{\pi}{4}$.
We can easily verify that $\mathbf{f}(t,\mathbf{x})$ in system (\ref{Corrspond-ODE2}) is continuous in $D$. Furthermore, $\mathbf{f}(t,\mathbf{x})$ is locally Lipschitzian with respect to $\mathbf{x}$ in $D$, which means for two arbitrary points
$(t,\mathbf{x}_1)=(t,\theta_1^1,\theta_2^1,\theta_3^1)$ and $(t,\mathbf{x}_2)=(t,\theta_1^2,\theta_2^2,\theta_3^2)$ in $D$, there is a $L$ such that $|\mathbf{f}(t,\mathbf{x}_2)-\mathbf{f}(t,\mathbf{x}_1)|<L|\mathbf{x}_2-\mathbf{x}_1|$.  It can be proved as follows. Since
\begin{equation*}
\begin{aligned}
&|\mathbf{f}(t,\mathbf{x}_2)-\mathbf{f}(t,\mathbf{x}_1)|
\\
\leq& |(\sin(2\theta_1^2) \tan(2\theta_2^2)-\sin(2\theta_1^1)\tan(2\theta_2^1)) a_2(t)-(\cos(2\theta_1^2)\tan(2\theta_2^2)-\cos(2\theta_1^1)\tan(2\theta_2^1))a_3(t)|
\\
&+|(\cos(2\theta_1^2)-\cos(2\theta_1^1))a_2(t)+(\sin(2\theta_1^2)-\sin(2\theta_1^1))a_3(t)|
\\
&+\left|\left(-\frac{\sin(2\theta_1^2)}{\cos (s2\theta_2^2)}+\frac{\sin(2\theta_1^1)}{\cos(2\theta_2^1)}\right)a_2(t)+\left(\frac{\cos(2\theta_1^2)}{\cos(2\theta_2^2)}-\frac{\cos(2\theta_1^1)}{\cos(2\theta_2^1)}\right)a_3(t)\right|
\\
\leq&|\sin(2\theta_1^2)\tan(2\theta_2^2)-\sin(2\theta_1^1)\tan(2\theta_2^1)||a_2(t)|+|\cos(2\theta_1^2)\tan(2\theta_2^2)-\cos(2\theta_1^1)\tan(2\theta_2^1)||a_3(t)|
\\
&+|\cos(2\theta_1^2)-\cos(2\theta_1^1)||a_2(t)|+|\sin(2\theta_1^2)-\sin(2\theta_1^1)||a_3(t)|
\\
&+\left|-\frac{\sin(2\theta_1^2)}{\cos(2\theta_2^2)}+\frac{\sin(2\theta_1^1)}{\cos(2\theta_2^1)}\right||a_2(t)|+\left|\frac{\cos(2\theta_1^2)}{\cos(2\theta_2^2)}-\frac{\cos(2\theta_1^1)}
{\cos(2\theta_2^1)}\right||a_3(t)|.
\end{aligned}
\end{equation*}
For the coefficients of the above inequality,
\begin{equation*}
\begin{aligned}
&|\sin(2\theta_1^2)\tan(2\theta_2^2)-\sin(2\theta_1^1)\tan(2\theta_2^1)|
\\
&=\left|\frac{\sin(2\theta_1^2)\sin(2\theta_2^2)\cos(2\theta_2^1)-\sin(2\theta_1^1)\sin(2\theta_2^1)\cos(2\theta_2^2)}{\cos(2\theta_2^2)\cos(2\theta_2^1)} \right|
\\
&=\frac{1}{|\cos(2\theta_2^1)\cos(2\theta_2^2)|}\left|\sin(2\theta_1^2)\sin(2\theta_2^2)\cos(2\theta_2^1)-\sin(2\theta_1^1)\sin(2\theta_2^1)\cos(2\theta_2^2)\right|
\\
&=\frac{1}{|\cos(2\theta_2^1)\cos(2\theta_2^2)|}\left|\sin(2\theta_1^2)\sin(2\theta_2^2)\cos(2\theta_2^1)-\sin(2\theta_1^1) \sin(2\theta_2^2) \cos(2\theta_2^1)\right.
\\
&~~+\sin(2\theta_1^1)\sin(2\theta_2^2)\cos(2\theta_2^1)-\sin(2\theta_1^1)\sin(2\theta_2^1)\cos(2\theta_2^1)
\\
&~~+\left.\sin(2\theta_1^1)\sin(2\theta_2^1)\cos(2\theta_2^1)-\sin(2\theta_1^1)\sin(2\theta_2^1)\cos(2\theta_2^2)\right|
\\
&\leq\frac{1}{|\cos(2\theta_2^1)\cos(2\theta_2^2)|}\left\{\left|\sin(2\theta_1^2)-\sin(2\theta_1^1)\right||\sin(2\theta_2^2)\cos(2\theta_2^1)|\right.
\\
&~~+\left|\sin(2\theta_2^2)-\sin(2\theta_2^1)\right||\sin(2\theta_1^1)\cos(2\theta_2^1)|
+\left.\left|\cos(2\theta_2^1)-\cos(2\theta_2^2)\right||\sin(2\theta_1^1)\sin(2\theta_2^1)|\right\}
\\
&\leq\frac{1}{|\cos(2\theta_2^1)\cos(2\theta_2^2|)}\big\{\left|2\cos(\theta_1^2+\theta_1^1)\sin(\theta_1^2-\theta_1^1)\right||\sin(2\theta_2^2)\cos(2\theta_2^1)|
\\
&~~+\left|2\cos(\theta_2^2+\theta_2^1)\sin(\theta_2^2-\theta_2^1)\right||\sin(2\theta_1^1)\cos(2\theta_2^1)|
\\
&~~
+\left|2\cos(\theta_2^1+\theta_2^2)\sin(\theta_2^1-\theta_2^2)\right||\sin(2\theta_1^1)\sin(2\theta_2^1)|\big\}
\\
&\leq\frac{1}{|\cos(2\theta_2^1)\cos(2\theta_2^2)|}\big\{|2\cos(\theta_1^2+\theta_1^1)\sin(2\theta_2^2)\cos(2\theta_2^1)|\big\}|\theta_1^2-\theta_1^1|
\\
&~~+\frac{1}{|\cos(2\theta_2^1)\cos(2\theta_2^2)|}\big\{\left|2\cos(\theta_2^2+\theta_2^1)\sin(2\theta_1^1)\cos(2\theta_2^1)\right|\big\}|\theta_2^2-\theta_2^1|
\\
&~~+\frac{1}{|\cos(2\theta_2^1)\cos(2\theta_2^2)|}\big\{|2\cos(\theta_2^2+\theta_2^1)\sin(2\theta_1^1)\sin(2\theta_2^1)|\big\}|\theta_2^2-\theta_2^1|.
\end{aligned}
\end{equation*}
According to the choice of $D$, all the triangle functions in the above inequalities
are continuous and bounded. So there is a real constant $k_1$ such that
\begin{equation*}
|\sin(2\theta_1^2)\tan(2\theta_2^2)-\sin(2\theta_1^1)\tan(2\theta_2^1)|\leq k_1 |\theta_1^2-\theta_1^1|+2k_1|\theta_2^2-\theta_2^1|.
\end{equation*}
For other coefficients, we similarly obtain
\begin{equation*}
\begin{aligned}
&|\cos(2\theta_1^2)\tan(2\theta_2^2)-\cos(2\theta_1^1)\tan(2\theta_2^1)|\leq k_2 |\theta_1^2-\theta_1^1|+2k_2|\theta_2^2-\theta_2^1|,
\\
&|\cos(2\theta_1^2)-\cos(2\theta_1^1)|\leq 2|\theta_1^2-\theta_1^1|,
\\
&|\sin(2\theta_1^2)-\sin(2\theta_1^1)|\leq 2|\theta_1^2-\theta_1^1|,
\\
&\left|-\frac{\sin(2\theta_1^2)}{\cos(2\theta_2^2)}+\frac{\sin(2\theta_1^1)}{\cos(2\theta_2^1)}\right|\leq k_3|\theta_1^2-\theta_1^1|+k_3|\theta_2^2-\theta_2^1|,
\\
&\left|\frac{\cos(2\theta_1^2)}{\cos(2\theta_2^2)}-\frac{\cos(2\theta_1^1)}{\cos(2\theta_2^1)}\right|\leq k_4|\theta_1^2-\theta_1^1|+k_4|\theta_2^2-\theta_2^1|,
\end{aligned}
\end{equation*}
where $k_2, k_3$ and $k_4$ are positive constants.
Since $a_1(t),a_2(t)$ and $a_3(t)$ are continuous on the closed set $D$, there exists a constant $M>0$ such that $|a_i(t)|\leq M$ for $i=1,2,3$. Hence
\begin{equation*}
\begin{aligned}
&|\mathbf{f}(t,\mathbf{x}_2)-\mathbf{f}(t,\mathbf{x}_1)|
\\
\leq&(k_1 |\theta_1^2-\theta_1^1|+2k_1|\theta_2^2-\theta_2^1|)M+(k_2 |\theta_1^2-\theta_1^1|+2k_2|\theta_2^2-\theta_2^1|)M
\\
+&(2|\theta_1^2-\theta_1^1|)M+(2|\theta_1^2-\theta_1^1|)M
\\
+&(k_3|\theta_1^2-\theta_1^1|+k_3|\theta_2^2-\theta_2^1|)M+(k_4|\theta_1^2-\theta_1^1|+k_4|\theta_2^2-\theta_2^1|)M
\\
\leq&(k_1+k_2+4+k_3+k_4)M|\theta_1^2-\theta_1^1|+(2k_1+2k_2+k_3+k_4)M|\theta_2^2-\theta_2^1|
\\
\leq&L\sqrt{(\theta_1^2-\theta_1^1)^2+(\theta_2^2-\theta_2^1)^2+(\theta_3^2-\theta_3^1)^2}=L|\mathbf{x}_2-\mathbf{x}_1|,
\end{aligned}
\end{equation*}
where $L:=(3k_1+3k_2+4+2k_3+2k_4)M$, that is,
$\mathbf{f}(t,\mathbf{x})$ is Lipschitzian with respect to $\mathbf{x}\in D$.
By Theorem 3.1 of chapter 1 in \cite{Hale1969}, the equation (\ref{Corrspond-ODE2})
with the initial value ${\bf x}(t_0)=(0,0,0)$
has a unique solution on the interval $I$
defined in the theorem.
\qquad$\Box$

\begin{remk}\label{singular-case-2}
{\rm
Similarly to Theorem 2.1  in \cite{Hale1969}, we have the following result of continuation:
If $D$ is an open set in $\mathbb{R}\times (-\pi,\pi)\times (-\frac{\pi}{4},\frac{\pi}{4})\times (-\frac{\pi}{2},\frac{\pi}{2})$
and $a_1(t),a_2(t)$, $a_3(t)$ are continuous functions in $D$,
then
the solution $\mathbf{x}$ of  the decisive equation (\ref{Corrspond-ODE2})
with the initial value ${\bf x}(t_0)=(0,0,0)$, obtained in Theorem \ref{existence},
can be extended
to a maximal interval
because $\mathbf{f}(t,\mathbf{x})$ is continuous in $D$.
Furthermore, if $(a,b)$ is a maximal interval of the solution, then $\mathbf{x}(t)$ tends to the boundary of $D$ as $t\rightarrow a$ and $t\rightarrow b$.
}
\end{remk}


\section{Computation in special cases}

Theorem \ref{existence} actually indicates an algorithm with Picard's iteration to compute
the solution of decisive equation (\ref{Corrspond-ODE2})
with the initial value ${\bf x}(t_0)=(0,0,0)$.
In this section, we give some results on exact solutions of
the decisive equation (\ref{Corrspond-ODE2}),
a nonautonomous nonlinear equation,
in the special cases:
{\bf (I)} $\theta_1(t)=0$,
{\bf (II)} $\theta_2(t)=0$,
{\bf (III)} $\theta_3(t)=0$.
For this purpose, we go back to QDE (\ref{PureQDE-IVP}) and find solutions of the forms
$q(t)=e^{\mathbf{j}\theta_2(t)}e^{\mathbf{k}\theta_3(t)}$,
$q(t)=e^{\mathbf{i}\theta_1(t)}e^{\mathbf{k}\theta_3(t)}$ and
$q(t)=e^{\mathbf{i}\theta_1(t)}e^{\mathbf{j}\theta_2(t)}$ respectively.


\subsection{Case (I)}

\begin{cor}\label{Corollary-1}
If
$a_1(t)=a_3(t)\tan(2A_2(t))$, where $A_2(t)=\int_{t_0}^t a_2(t)dt$, QDE (\ref{PureQDE-IVP}) has a solution $q(t)=e^{\mathbf{j}A_2(t)}e^{\mathbf{k}\theta_3}\mathbf{C}_0$, where
$\theta_3(t):=\int_{t_0}^t\frac{a_3(t)}{\cos(2A_2(t))}dt$ and $\mathbf{C}_0$ is a constant.
\end{cor}

{\bf Proof.}  For the special case $\theta_1(t)=0$, i.e., $q(t)=e^{\mathbf{j}\theta_2}e^{\mathbf{k}\theta_3}$,
the corresponding system (\ref{Corrspond-ODE2}) is
\begin{equation*}
\left\{\begin{aligned}
\theta'_1(t)&= a_1(t)-\tan(2\theta_2)a_3(t)=0 ,
\\
\theta'_2(t)&=a_2(t),
\\
\theta'_3(t)&=\frac{1}{\cos(2\theta_2)}a_3(t),
\end{aligned}\right.
\end{equation*}
which means $\theta_2(t)=\int_{t_0}^t a_2(t)dt=A_2(t)$ and
$a_1(t)=a_3(t)\tan(2A_2(t))$. Under this condition, one can check that
the function $(\theta_1(t),\theta_2(t),\theta_3(t))$, where
\begin{eqnarray*}
\left\{
\begin{aligned}
\theta_1(t)&=0,
\\
\theta_2(t)&=A_2(t),
\\
\theta_3(t)&=\int_{t_0}^t\frac{a_3(t)}{\cos(2\theta_2(t))}dt=\int_{t_0}^t\frac{a_3(t)}{\cos(2A_2(t))}dt,
\end{aligned}
\right.
\end{eqnarray*}
is a special solution of (\ref{Corrspond-ODE2}). Hence  the general solution of  QDE (\ref{PureQDE-IVP}) is $q(t)=e^{\mathbf{j}\theta_2}e^{\mathbf{k}\theta_3}\mathbf{C}_0$ and the proof is completed.
 \qquad$\Box$


\subsection{Case (II)}

\begin{cor}\label{Corollary-2}
If $a_2(t)=-a_3(t)\tan(2A_1(t))$, where $A_1(t)=\int_{t_0}^t a_1(t)dt$, QDE (\ref{PureQDE-IVP}) has a solution $q(t)=e^{\mathbf{i}A_1(t)}e^{\mathbf{k}\theta_3(t)}\mathbf{C}_0$, where
$\theta_3(t):=\int_{t_0}^t\frac{a_3(t)}{\cos(2A_1(t))}dt$ and $\mathbf{C}_0$ is a constant.
\end{cor}

{\bf Proof.}
For the special case $\theta_2(t)=0$, i.e., $q(t)=e^{\mathbf{i}\theta_1}e^{\mathbf{k}\theta_3}$,
the corresponding system (\ref{Corrspond-ODE2}) is
\begin{equation*}
\left\{\begin{aligned}
\theta'_1(t)&= a_1(t) ,
\\
\theta'_2(t)&=\cos(2\theta_1)a_2(t)+\sin(2\theta_1)a_3(t)=0,
\\
\theta'_3(t)&=-\sin(2\theta_1)a_2(t)+\cos(2\theta_1)a_3(t),
\end{aligned}\right.
\end{equation*}
which means $\theta_1(t)=\int_{t_0}^t a_1(t)dt=A_1(t)$ and
$a_2(t)=-a_3(t)\tan(2A_1(t))$.
Under this conditions, one can check that
the function $(\theta_1(t),\theta_2(t),\theta_3(t))$, where
\begin{eqnarray*}
\left\{
\begin{aligned}
\theta_1(t)&=A_1(t),
\\
\theta_2(t)&=0,
\\
\theta_3(t)&=\int_{t_0}^t\frac{a_3(t)}{\cos(2A_1(t))}dt,
\end{aligned}
\right.
\end{eqnarray*}
 is a solution of ODEs (\ref{Corrspond-ODE2}). Hence QDE (\ref{PureQDE-IVP}) has a solution $q(t)=e^{\mathbf{i}\theta_1(t)}e^{\mathbf{k}\theta_3(t)}\mathbf{C}_0$
 and the proof is completed.
\qquad$\Box$


\subsection{Case (III)}
\begin{cor}\label{Corollary-3}
If $a_3(t)=a_2(t)\tan2A_1(t)$, where $A_1(t)=\int_{t_0}^t a_1(t)dt$, QDE (\ref{PureQDE-IVP}) has a solution $q(t)=e^{\mathbf{i}A_1(t)}e^{\mathbf{j}\theta_2(t)}\mathbf{C}_0$, where $\theta_2(t):=\int_{t_0}^t \frac{a_2(t)}{\cos (2A_1(t))}dt$ and $\mathbf{C}_0$ is a constant.
\end{cor}

{\bf Proof.} For the special case $\theta_3(t)=0$, i.e., $q(t)=e^{\mathbf{i}\theta_1}e^{\mathbf{j}\theta_2}$,
the corresponding system (\ref{Corrspond-ODE2}) is
\begin{equation*}
\left\{\begin{aligned}
\theta'_1(t)&= a_1(t)+\sin(2\theta_1)\tan(2\theta_2)a_2(t)-\cos(2\theta_1)\tan(2\theta_2) a_3(t) ,
\\
\theta'_2(t)&=\cos(2\theta_1)a_2(t)+\sin(2\theta_1)a_3(t),
\\
\theta'_3(t)&= -\frac{\sin(2\theta_1)}{\cos (2\theta_2)}a_2(t)+\frac{\cos(2\theta_1)}{\cos(2\theta_2)}a_3(t)=0,
\end{aligned}\right.
\end{equation*}
which means $\theta_1(t)=\int_{t_0}^t a_1(t)dt=A_1(t)$ and $a_3(t)=a_2(t)\tan(2A_1(t))$.
Under this conditions, one can check that
the function $(\theta_1(t),\theta_2(t),\theta_3(t))$, where
\begin{eqnarray*}
\left\{
\begin{aligned}
\theta_1(t)&=A_1(t),
\\
\theta_2(t)&=\int_{t_0}^t \frac{a_2(t)}{\cos (2A_1(t))}dt,
\\
\theta_3(t)&=0,
\end{aligned}
\right.
\end{eqnarray*}
 is a solution of ODEs (\ref{Corrspond-ODE2}). Hence QDE (\ref{PureQDE-IVP}) has a solution $q(t)=e^{\mathbf{i}\theta_1(t)}e^{\mathbf{j}\theta_2(t)}\mathbf{C}_0$
 and the proof is completed.
\qquad$\Box$

\section{Examples}

\begin{exam}\label{Ex1}
Solve QDE $q'(t)=(\sin (2t) \mathbf{i}+\mathbf{j}+\cos(2t)\mathbf{k})q(t)$.
\end{exam}

Note that
$a(t):=\sin (2t) \mathbf{i}+\mathbf{j}+\cos(2t)\mathbf{k}=\underline{a}(t)$, which
does not satisfy condition (\ref{aaabbb}). So  it cannot be solved by the formula
of QDE with commutativity.

Since
$a(t)=a_0(t)+a_1(t) {\bf i} +a_2(t){\bf j}+a_3(t){\bf k}$ as denoted before, we have
$a_1(t)=\sin (2t), a_2(t)=1$ and $a_3(t)=\cos (2t)$.
Then, from (\ref{a-A}) we obtain
$A_2(t)=\int_0^t a_2(t)dt=t$ and $a_1(t)=\sin(2A_2(t))=\cos(2A_2(t))\tan(2A_2(t))=a_3\tan(2A_2(t)).$
By Corollary \ref{Corollary-1}, the decisive equation (\ref{Corrspond-ODE2}) has
a solution $(\theta_1,\theta_2,\theta_3)^T$, where
\begin{equation*}
\left\{
\begin{aligned}
\theta_1(t)&=0,
\\
\theta_2(t)&=t,
\\
\theta_3(t)&=\int_0^t g(t)dt=t.
\end{aligned}
\right.
\end{equation*}
This implies that the general solution of the QDE is
$q(t)=e^{\mathbf{j}t}e^{\mathbf{k}t}\mathbf{C}$
, where $\mathbf{C}$ is an arbitrary quaternionic number.

Further, we compute
\begin{equation*}
\begin{aligned}
q'(t)&=\mathbf{j}e^{\mathbf{j}t}e^{\mathbf{k}t}+e^{\mathbf{j}t}\mathbf{k}e^{\mathbf{k}t}
\\
&=\mathbf{j} \left(\cos (t)+\mathbf{j}\sin (t)\right) \left(\cos (t)+\mathbf{k}\sin (t)\right)+ \left(\cos (t)+\mathbf{j}\sin (t)\right)\mathbf{k} \left(\cos (t)+\mathbf{k}\sin (t)\right)
\\
&=\left(-\cos (t) \sin (t)+\mathbf{i}\cos (t) \sin (t)+\mathbf{j}\cos^2(t)-\mathbf{k}\sin^2 (t) \right)
\\
&+\left(-\cos (t) \sin (t)+\mathbf{i}\cos (t) \sin (t)-\mathbf{j}\sin^2(t)+\mathbf{k}\cos^2 (t)\right)
\\
&=\left(-\sin(2t) +\mathbf{i}\sin(2t)+\mathbf{j}\cos(2t)+\mathbf{k}\cos(2t)\right)
\end{aligned}
\end{equation*}
and
\begin{equation*}
\begin{aligned}
a(t)q(t)&=\left(\sin (2t) \mathbf{i}+\mathbf{j}+\cos (2t) \mathbf{k} \right)q(t)= \left(\sin (2t) \mathbf{i}+\mathbf{j}+\cos (2t) \mathbf{k}\right)e^{\mathbf{j}t}e^{\mathbf{k}t}
\\
&= \left(\sin(2t) \mathbf{i}+\mathbf{j}+\cos(2t) \mathbf{k}\right) \left(\cos (t)+\mathbf{j}\sin (t)\right) \left(\cos (t)+\mathbf{k}\sin (t)\right)
\\
&=\left(\sin (2t) \mathbf{i}+\mathbf{j}+\cos (2t)\mathbf{k} \right) \left(\cos^2(t)+\mathbf{i}\sin^2 (t)+\mathbf{j}\cos (t) \sin (t)+\mathbf{k}\cos (t) \sin (t)\right)
\\
&=\left(-\sin (2t) \sin^2(t)-\sin (t) \cos (t)-\cos(2t) \sin (t) \cos (t)\right)
\\
&+\mathbf{i} \left(\sin (2t) \cos^2(t)+\sin (t) \cos (t)-\cos(2t)\sin (t)\cos (t)\right)
\\
&+\mathbf{j} \left(-\sin (2t) \sin (t) \cos (t)+\cos^2 (t)+\cos (2t) \sin^2(t)\right)
\\
&+\mathbf{k}(\sin (2t) \sin (t) \cos (t) -\sin^2(t)+\cos(2t)\cos^2(t))
\\
&=\left(-2\sin^3(t)\cos (t)-\sin (t)\cos (t)-2\sin (t) \cos^3(t)+\sin (t)\cos (t)\right)
\\
&+\mathbf{i}\left(2\sin (t) \cos^3(t)+\sin (t)\cos (t)-2\sin (t) \cos^3 (t)+\sin (t) \cos (t)\right)
\\
&+\mathbf{j} \left(-2\sin^2(t)\cos^2(t)+\cos^2(t)+\sin^2(t) \cos^2(t)-\sin^4(t)\right)
\\
&+\mathbf{k}(2\sin^2(t) \cos^2(t)-\sin^2(t)+\cos^4(t)-\sin^2(t)\cos^2(t))
\\
&=(-\sin (2t)+\mathbf{i}\sin (2t)+\mathbf{j}\cos(2t)+\mathbf{k}\cos(2t)).
\end{aligned}
\end{equation*}
The same result of the above computation verifies that $q$ satisfies the QED.

\begin{exam}Solve QDE $q'(t)=(\mathbf{i}+t\sin (2t) \mathbf{j}-t\cos (2t) \mathbf{k})q(t)$.
\end{exam}




Note that $a(t)=\underline{a}(t)=\mathbf{i}+t\sin (2t) \mathbf{j}-t\cos (2t) \mathbf{k}$, which does not satisfy condition
(\ref{aaabbb}),
i.e., the commutativity condition is not true.
So this equation cannot be solved by the formula
of QDE with commutativity.

Since
$a(t)=\underline{a}(t)=a_1(t) {\bf i} +a_2(t){\bf j}+a_3(t){\bf k}$ as denoted before,
we have $a_1(t)=1$, $ a_2(t)=t\sin(2t)=t\sin(2A_1) $ and $a_3(t)=-t\cos(2A_1)$, where $A_1=\int_0^t a_1(t)dt=t$.
Hence $ a_2(t)=-a_3(t)\tan(2A_1)$.
By Corollary \ref{Corollary-2}, the solution of system (\ref{Corrspond-ODE2}) is $(\theta_1,\theta_2,\theta_3)^T$, where
\begin{equation*}
\left\{
\begin{aligned}
\theta_1(t)&=t,
\\
\theta_2(t)&=0,
\\
\theta_3(t)&=\int_{0}^t\frac{a_3(t)}{\cos (2A_1(t))}dt=-\frac{t^2}{2}.
\end{aligned}
\right.
\end{equation*}
Thus
the QDE has the solution $q(t)=e^{\mathbf{i}t}e^{-\mathbf{k}\frac{t^2}{2}}\mathbf{C}$, where $\mathbf{C}$ is an arbitrary quaternionic number.
One can compute
\begin{equation*}
\begin{aligned}
q'(t)
&=\left[-\sin (t) \cos \left( \frac{t^2}{2}\right) -t\cos (t) \sin \left(\frac{t^2}{2} \right) \right]+\left[\cos (t)\cos \left(\frac{t^2}{2}\right)-t\sin (t) \sin \left( \frac{t^2}{2}\right) \right]\mathbf{i}
\\
&+\left[\cos (t) \sin \left(\frac{t^2}{2}\right)+t\sin (t) \cos \left(\frac{t^2}{2}\right) \right]\mathbf{j}+\left[\sin (t) \sin\left( \frac{t^2}{2}\right)-t\cos (t) \cos \left(\frac{t^2}{2}\right)\right]\mathbf{k}
\\
&=(\mathbf{i}+t \sin (2t) \mathbf{j}-t\cos (2t) \mathbf{k})q(t),
\end{aligned}
\end{equation*}
which checks that this function $q(t)$ satisfies the QED.

\begin{exam}Solve QDE $q'(t)=(\mathbf{i}+t\cos (2t)\mathbf{j}+t\sin(2t)\mathbf{k})q(t)$.
\end{exam}

Note that $a(t)=\underline{a}(t)=\mathbf{i}+t\cos (2t) \mathbf{j}+t\sin(2t)\mathbf{k}$, which does not satisfy condition
(\ref{aaabbb}),
i.e., the commutativity condition is not true.
So this equation cannot be solved by the formula
of QDE with commutativity.

Since
$a(t)=\underline{a}(t)=a_1(t) {\bf i} +a_2(t){\bf j}+a_3(t){\bf k}$ as denoted before,
we have $a_1(t)=1$, $ a_2(t)=t\cos(2t)=t\cos(2A_1) $ and $a_3(t)=t\sin(2A_1)$, where $A_1=\int_0^t a_1(t)dt=t$.
Hence $a_3(t)=a_2(t)\tan(2A_1(t))$.
By Corollary \ref{Corollary-3}, the solution of ODEs (\ref{Corrspond-ODE2}) is $(\theta_1,\theta_2,\theta_3)^T$, where
\begin{equation*}
\left\{
\begin{aligned}
\theta_1(t)&=t,
\\
\theta_2(t)&=\int_{0}^t \frac{a_2(t)}{\cos(2A_1(t))}dt=\frac{t^2}{2},
\\
\theta_3(t)&=0.
\end{aligned}
\right.
\end{equation*}
Thus the QDE has the solution $q(t)=e^{\mathbf{i}t}e^{\mathbf{j}\frac{t^2}{2}}\mathbf{C}$, where $\mathbf{C}$ is an arbitrary quaternionic number.
One can compute
\begin{equation*}
\begin{aligned}
q'(t)
&=\left[-\sin (t)\cos \left(\frac{t^2}{2}\right)-t\cos( t)\sin \left(\frac{t^2}{2}\right)\right]+\left[\cos (t)\cos \left(\frac{t^2}{2}\right)-t\sin (t) \sin\left(\frac{t^2}{2}\right)\right]\mathbf{i}
\\
&+\left[-\sin (t) \sin \left(\frac{t^2}{2}\right)+t\cos (t) \cos\left(\frac{t^2}{2}\right)\right]\mathbf{j}+\left[\cos (t)\sin \left(\frac{t^2}{2}\right)+t\sin (t)\cos\left(\frac{t^2}{2}\right)\right]\mathbf{k}
\\
&=(\mathbf{i}+t\cos (2t)\mathbf{j}+t\sin(2t)\mathbf{k})q(t),
\end{aligned}
\end{equation*}
which checks that this function $q(t)$ satisfies the QED.

\vskip 0.2cm

{\bf Acknowledgement:}
This paper started in an collaboration in the University of Macau.
All authors are ranked in alphabetic order of their names and their contributions should be treated equally.

\end{document}